\documentclass[10pt]{amsart}
\usepackage{amsmath, amsfonts, amsthm, amssymb, verbatim,dsfont}
\usepackage{graphicx}
\usepackage[mathcal]{euscript}
\usepackage{amstext}
\usepackage{dsfont} 
\usepackage{amssymb,mathrsfs}
\usepackage[usenames]{color}
\theoremstyle{plain}
        \newtheorem{theorem}{Theorem}[section]
        \newtheorem{proposition}[theorem]{Proposition}
        \newtheorem{lemma}[theorem]{Lemma}
        
\theoremstyle{definition}

\theoremstyle{plain}
        
\numberwithin{equation}{section}

\usepackage{amsmath}
\usepackage{verbatim}
\newcommand \be           {\begin{equation}}
\newcommand \ee            {\end{equation}}
\newcommand \Dcal           {\mathcal{D}}
\newcommand \Ecal           {\mathcal{E}}

\newcommand \Mcal           {\mathcal{M}}
\newcommand \Qcal           {\mathcal{Q}}
\newcommand \RR           {\mathbb{R}}
\newcommand \NN           {\mathbb{N}}

\newcommand \CC           {\mathbb{C}}

\newcommand \Pbold           {\mathbf{P}} 
\newcommand \PP \Pbold

\newcommand \del           \partial
\newcommand \eps            \epsilon

\newcommand \OO     {\mathcal{O}}

\newcommand \loc        {{\mathrm{loc}}}

\DeclareMathOperator    \Real {Re}
\DeclareMathOperator    \Imag {Im}
\newcommand \ws     {\mathrel{\mathop{\rightharpoonup}\limits^{*}}}
\definecolor{gray}{gray}{0.4}
\newcommand{\unj}{u^n_j}

\newcommand{\unnj}{u^{n+1}_j}

\newcommand{\vnj}{v^n_j}

\newcommand{\vnnj}{v^{n+1}_j}

\newcommand{\vh}{v^h}
\newcommand{\uh}{u^h}
\newcommand{\vj}{v_j}

\newcommand{\zj}{z_j}
\newcommand{\zjm}{z_{j-1}}
\newcommand{\zjmm}{z_{j-2}}
\newcommand{\zjj}{z_{j+1}}
\newcommand{\zjjj}{z_{j+2}}
\newcommand{\uj}{u_j}

\newcommand{\ujj}{u_{j+1}}

\newcommand{\xj}{x_j}

\newcommand{\xjj}{x_{j+1}}

\def\XXint#1#2#3{{\setbox0=\hbox{$#1{#2#3}{\int}$}
\vcenter{\hbox{$#2#3$}}\kern-.5\wd0}}

\def\build#1_#2^#3{\mathrel{
\mathop{\kern 0pt#1}\limits_{#2}^{#3}}}
\begin{document}

\title[a numerical scheme for a coupled Sch\-r\"o\-din\-ger--KdV]{Convergence of a numerical scheme for a coupled Sch\-r\"o\-din\-ger--KdV system}

\author[P. Amorim, M. Figueira]{Paulo Amorim$^1$, M\'ario Figueira$^1$}

\date{\today}


\footnotetext[1]{Centro de Matem\'atica e Aplica\c c\~oes
Fundamentais, Universidade de Lisboa, Av. Prof. Gama Pinto 2,
1649-003 Lisboa, Portugal. E-mail: {\tt pamorim@ptmat.fc.ul.pt, figueira@ptmat.fc.ul.pt}}

\maketitle

\begin{abstract}
We prove the convergence in a strong norm 
of a finite difference semi-discrete scheme approximating a coupled
Sch\-r\"o\-din\-ger--KdV system on a bounded domain. This system models the interaction of short and long waves. 
Since the energy estimates available in the continuous case do not carry over to the
discrete setting, we rely on a suitably truncated problem which we prove reduces to the
original one. We present some numerical examples to illustrate our convergence result.
\\[5pt]
\end{abstract}

\section{Introduction}
In \cite{Benney}, D.J. Benney presents a general theory modeling the nonlinear interaction between 
short waves and long waves, deriving 
nonlinear differential systems describing these interactions in various physical settings.
The (complex-valued) short waves $u(x,t)$, $x\in\RR$, $t\ge0$, are
described by a nonlinear Sch\-r\"o\-din\-ger equation and the (real-valued) long waves $v(x,t)$ satisfy a 
quasilinear equation, eventually with a dispersive term. 
In the most general context, the interaction is described by the nonlinear system
\[
\left\{
\aligned 
& i \del_t u + i c_1 \del_x u + \del_{xx} u= \alpha vu + \beta |u|^2 u
\\
& \del_t v +  c_2 \del_x v + \mu \del^3_x v + \nu \del_x v^2 = \gamma \del_x (|u|^2),
\endaligned \right.
\]
where $c_1,c_2, \alpha, \beta, \gamma, \mu$ and $\nu$ are real constants.

In this paper, we are concerned with the numerical approximation of the solutions to the Cauchy problem 
for the system comprising the nonlinear Sch\-r\"o\-din\-ger equation coupled with a 
Korteweg--de Vries equation with Dirichlet boundary conditions on a bounded domain $(0,L)$,
\begin{subequations}
\label{15}
\begin{align}
\label{20}
&i\del_t u + \del_{xx} u = \alpha v u + \beta |u|^2 u 
\\
\label{30}
& \del_t v +  \del_{x}^3 v +  \del_x(v^2) = \gamma \del_x| u|^2
\\
\label{40}
& u(x,0) = u_0(x), \quad v(x,0) = v_0(x), \qquad x\in (0,L), \quad L>0,
\\ & u(0,t) = u(L,t) = 0, \quad v(0,t) = v(L,t) = 0, \quad t\in [0,T], \quad T>0.
\end{align}
\end{subequations}

This kind of system arises in several fields of physics such as the study of resonant 
interactions, short and long capillary-gravity waves on water \cite{8}, an electron plasma 
interaction \cite{9} and a diatomic lattice system \cite{10}.

The Cauchy problem for the system \eqref{15} was initially studied on the whole line by M.~Tsutsumi
\cite{11}, who proved that for initial data $(u_0, v_0)$ in $H^{m+1/2}(\RR)\times H^m (\RR)$, $m=1,2.\dots,$ 
the problem is well-posed in the same space. After that, Bekiranov \emph{et al.~}\cite{12} established local well-posedness for initial data in $H^s(\RR)\times H^{s-1/2}(\RR), s\ge 0$, and more recently, Corcho
and Linares \cite{13} proved global well-posedness in the energy space $H^1(\RR) \times H^1(\RR)$.

It is worth pointing out that in the above well-posedness results, uniqueness is obtained only in 
some subspace of $C([0,T]; H^s(\RR)\times H^{s-1/2}(\RR))$ (in \cite{12}), (resp.~a subspace of
$C([0,T]; H^1(\RR)\times H^1(\RR))$ in \cite{13}). Additionally, the techniques used in the 
papers\cite{12,13} (which were introduced by Bourgain \cite{14} and developed by Kenig, Ponce 
and Vega \cite{15,16}),
do not seem to be applicable to the generalized KdV equation (gKdV). In a very recent paper, Dias \emph{et al.~}\cite{17}, using energy methods, obtain a global solution in $(H^1(\RR))^2$
for a coupled Sch\-r\"o\-din\-ger--gKdV system.

In this paper, we prove a convergence result for a semi-discrete finite difference approximation
of the system \eqref{15} in the space $(H^1(0,L))^2$.
The energy methods used by M.~Tsutsumi \cite{11} to prove 
global existence of a solution fail in the discrete setting, so we propose a new approach: 
by an appropriate truncation of the quadratic function $v^2$ appearing in the equation \eqref{30}, 
we consider a related problem involving a gKdV equation. For each fixed value of the 
truncation parameter, 
we are able to prove the convergence of a numerical scheme toward the solutions of this
auxiliary problem. These solutions, in turn, satisfy an energy inequality. Lastly, using this energy inequality,
we are able to derive an $L^\infty$ estimate independent of the truncation parameter,
which implies that the truncated problem 
in fact reduces to the original one. 

In contrast to previous work, the proof of these stability estimates 
require that we work on a bounded subset $(0,L)\subset \RR$. From the viewpoint of the applications
and numerical approximation,
this presents no great loss in generality.

Note also that our convergence proof does not rely on any previous existence results, and so 
constitutes a new existence proof for the Cauchy problem \eqref{15}. Additionally, the present work represents,
to the authors' best knowledge, the first numerical treatment of the system \eqref{15}--\eqref{40}.

An outline of the paper follows. After some notations and preliminaries, we state in Section~\ref{110} 
the main convergence result, Theorem~\ref{225}. In Section~\ref{235}, we prove 
Theorem~\ref{225} and present Proposition~\ref{305}, our main auxiliary result, dealing with 
the convergence of approximate solutions to a suitably truncated system, and an energy estimate.
Its proof is the object of Section~\ref{335}. Finally,
in the last section of the paper we illustrate our convergence result with some numerical simulations
and check its accuracy by testing it against some known exact solutions.

\subsection{Notations and preliminaries}
Let us introduce the Banach spaces
\[
\aligned
X_{J,\CC} = \{ z^h = (z_0, z_1, \dots, z_{J+1}) \in \CC^{J+2} : z_0 = z_1 = z_J = z_{J+1} = 0\}
\endaligned
\]
with $J\in \NN_0$ and $h = L/(J+1)$. In a similar way, we define the real space $X_{J,\RR}$. 
When no ambiguity arises, we represent either of these spaces by $X_J$.
The scalar product is given by
\[
\aligned
(z^h, w^h) = \sum_{j=2}^{J-1} h \zj \overline w_j, \qquad z^h, w^h \in X_J,
\endaligned
\]
and the $p$-norms by
\[
\aligned
\| z^h \|_{p,h} = \Big( \sum_{j=2}^{J-1} h| z_j|^p \Big)^{1/p},\quad 1\le p <\infty;\quad
\| z^h \|_\infty = \max_{j=2,\dots,J-1} |z_j|, \quad z^h \in X_J.
\endaligned
\]
To simplify the notation we write $\| z \|_p$ for the norm of $z$ in both the continuous and 
the discrete case. We denote by $H^m(0,L),$ $H^m_0(0,L)$ and $H^{-m}(0,L)$ $(m\in \NN)$ 
the usual Sobolev spaces. All the norms appearing in this paper are in $(0,L)$, so for instance
$\| u\|_2$ means $\| u \|_{L^2(0,L)}$.

We will use the following notations for the finite difference operators. For $z = (z_j)$,
\[
D_+ z_j = \frac{\zjj - \zj}h, \qquad D_- \zj = \frac{\zj - \zjm} h,
\]
\[
D_0 \zj = \frac{\zjj - \zjm}{2h} = \frac12 (D_+ + D_-) \zj,
\]
\[
\Delta^h \zj =  D_- D_+ \zj  = D_+ D_- \zj = \frac{\zjj -2 \zj + \zjm}{h^2}, 
\]
\[
D^3 \zj = D_0 D_-D_+ \zj = \frac{\zjjj -2 \zjj +2 \zjm - \zjmm}{2h^3}.
\]
For $u = (\uj)$, let us now introduce the piecewise linear interpolator,
\be
\label{50}
\aligned
\PP_1^h u (x) = \uj + (x -x_j) \frac{\ujj - \uj}{\xjj - \xj}, \quad x\in (\xj, \xjj),
\endaligned
\ee
and the piecewise constant interpolator,
\[
\aligned
\PP_0^h u (x) = \uj,\quad x\in (\xj, \xjj).
\endaligned
\]
Let $\uh \in X_J$. From \eqref{50} we derive
\be
\label{60}
\aligned
\| \PP_1^h \uh - \PP_0^h \uh\|_2 \le C h\| D_+ \uh\|_2
\endaligned
\ee
for some $C$ independent of $h$. As a consequence, we obtain
\be
\label{70}
\aligned
\| \PP_1^h\uh \|_2 \le \| \PP_1^h \uh - \PP_0^h \uh\|_2 + \| \PP_0^h \uh\|_2 \le
C \|\uh\|_2.
\endaligned
\ee

The following lemma establishes some inequalities which will be of use throughout.
\begin{lemma}
\label{75}
Let $\phi = (\phi_j) \in X_J$. Then
\be
\label{80}
\|\phi \|_\infty \le \sqrt{2} \|\phi\|_2^{1/2} \|D_\pm \phi\|_2^{1/2},
\ee
\be
\label{90}
\|\phi \|_\infty \le \frac12 \|D_\pm \phi\|_1 \le \frac{\sqrt{L}}2 \| D_\pm \phi \|_2
\ee
\be
\label{100}
\|\phi \|_2 \le C \|\phi\|_\infty,
\ee
where $C$ is a constant independent of $h$.
\end{lemma}
\proof The inequality \eqref{80} is a consequence of the Gagliardo--Nirenberg inequality and 
$\del_x \PP_1^h = \PP_0^h D_+$:
\[
\aligned
\| \phi \|_\infty = \| \PP_1^h \phi \|_\infty &\le \sqrt{2} \| \PP_1^h \phi \|_2^{1/2} 
\| \del_x \PP_1^h \phi \|_2^{1/2}
\le  \sqrt{2}\|\phi \|_2^{1/2} \| D_\pm\phi\|_2^{1/2},
\endaligned
\]
while \eqref{90} is a consequence of the (continuous) inequality $\| \phi\|_\infty 
\le \frac12 \|\phi'\|_1$. The last inequality is an elementary consequence of the definition
of the discrete norms.
\endproof

\section{Statement of the main result}
\label{110}
We propose the following semidiscrete finite difference approximation to the 
Cauchy problem \eqref{15}:
\begin{subequations}
\label{115}
\begin{align}
\label{120}
&i\del_t \uh + \Delta^h \uh = \beta |\uh|^2 \uh + \alpha\vh \uh
\\
\label{130}
& \del_t \vh +  D^3\vh + D_0 ({\vh})^2 = \gamma D_0 | \uh|^2
\\
\label{140}
& \uh(0) = \uh_0, \quad \vh(0) = \vh_0,
\\
\label{145}
& \uh(t), \vh(t) \in X_J, \quad t\in [0,T].
\end{align}
\end{subequations}

The 
global existence proof of Tsutsumi \cite{11} relies on energy methods which we
cannot carry over to 
the finite difference framework. It turns out that in the semidiscrete case, the crux of our convergence
argument relies on an \emph{a priori} $L^\infty$ bound. However, this bound is only available for a 
modified problem (see Proposition \ref{305} below). To 
deal with this difficulty, 
we use the fact that, under the right conditions, this problem reduces to the 
original one.

Our main result establishes the convergence of the approximations \eqref{115} toward a 
global weak solution of the problem \eqref{15} in the space $(H^1(0,L))^2$.
\begin{theorem}
\label{225}
Let $\alpha, \beta, \gamma$ be such that $\alpha \gamma>0$, $(u_0, v_0) \in (H^1_0(0,L))^2$ and
let $(\uh, \vh) \in (C([0,T]; X_J))^2$ be the global solutions of the discretized problem \eqref{115}
with initial data $(u^h_0, v^h_0)$, such that $ \PP_1^h u^h_0 \to u_0$ and $ \PP_1^h v^h_0 \to v_0$ 
in $H^1(0,L)$ as $h\to 0$.
Then, up to a subsequence,
\[
\PP_1^h \uh \ws u, \quad \PP_1^h \vh \ws v \qquad \text{ in }\quad
L^\infty\big( [0,T] ; H^1(0,L) \big),
\]
\[
\PP_1^h \uh \to u, \quad \PP_1^h \vh \to v \qquad \text{ in }\quad
L^\infty\big( [0,T] ; L^2(0,L) \big),
\]
with 
\[
\aligned
(u, v) \in {}&\big( L^\infty ([0,T]; H^1_0(0,L)) \big)^2 
\\
&\cap \big( C( [0,T] ; H^1_0(0,L)) \times 
C( [0,T] ; L^2(0,L)) \big), \quad T>0,
\endaligned
\]
a weak solution of the 
Sch\-r\"o\-din\-ger--KdV system \eqref{15}.
\end{theorem}

\section{Proof of Theorem \ref{225}  and stability estimates}
\label{235}
In this section we prove Theorem~\ref{225} along with the necessary stability estimates. 
We begin with the definition of an appropriate truncated
problem.
For each $M>1$ we define $C^\infty$ functions $f^M, g^M$ satisfying
\[
f^M (v) = \left\{
\aligned
&v^2,  &&\text{ if } |v| \le M,
\\
& |v|,  &&\text{ if } |v| > M^2 + 1,
\endaligned
\right.
\]
and 
\[
g^M (v) = \left\{
\aligned
&v,  &&\text{ if } |v| \le M,
\\
& \pm C,  &&\text{ if } |v| > 2M,
\endaligned
\right.
\]
with $0\le f^M(v) \le v^2$.
Here, the constant $C=C(M)$ is chosen to ensure the following property,
\be
\label{150}
|(f^M)'|_\infty + |g^M|_\infty \le C(M), \quad |(g^M)'|_\infty \le 1.
\ee
The functions $f^M, g^M$ are simply appropriate smooth truncations of the functions $v^2$ and $v$
appearing in \eqref{30}. We define also $F^M(v) = \int_0^v f^M(s) \,ds.$ 

Now, we consider the auxiliary Cauchy problem,
\begin{subequations}
\label{155}
\begin{align}
\label{160}
&i\del_t u + \del_{xx} u =\beta |u|^2 u + \alpha g^M(v) u
\\
\label{170}
& \del_t v +  \del^3_{x} v +  \del_x f^M(v) = \gamma \del_x( (g^M)'(v)| u|^2)
\\
\label{180}
& u(x,0) = u_0(x), \quad v(x,0) = v_0(x), \qquad x\in (0,L), \quad L>0,
\\ & u(0,t) = u(L,t) = 0, \quad v(0,t) = v(L,t) = 0, \quad t\in [0,T], \quad T>0.
\end{align}
\end{subequations}
and we propose the following semidiscrete finite difference approximation of \eqref{155}:
\begin{subequations}
\label{185}
\begin{align}
\label{190}
&i\del_t \uh + \Delta^h \uh = \beta |\uh|^2 \uh +\alpha g^M(\vh) \uh
\\
\label{200}
& \del_t \vh +  D^3\vh + D_0 f^M ({\vh}) = \gamma D_0( (g^M)'(\vh) | \uh|^2)
\\
\label{210}
& \uh(0) = \uh_0, \quad \vh(0) = \vh_0,
\\
\label{220}
& \uh(t), \vh(t) \in X_J, \quad t\in [0,T].
\end{align}
\end{subequations}


We will need some conservation laws of the auxiliary system 
\eqref{155}. A first result is as follows:
\begin{proposition}
\label{245}
Let $I\subset \RR$ be an interval and let $(u,v) \in (L^\infty_\loc (I ; H^1_0(0,L)))^2 $  be a solution
of the auxiliary system \eqref{155}. Then, for each $t,s \in I$, we have 
\be
\label{250}
\Mcal^M (t)  :=  \| u(t) \|_2 = \| u(s) \|_2 = \Mcal^M(s)
\ee
\be
\label{260}
\Qcal^M (t)  :=  \alpha \int_0^L v^2(x,t) \,dx + 2 \gamma \Imag \int_0^L u(x,t) 
\del_x \overline u(x,t) \,dx 
= \Qcal^M(s).
\ee
\end{proposition}
\proof 
Although the estimates \eqref{250},\eqref{260} are formally easy to obtain, the rigorous justification
of \eqref{260} is nontrivial and requires techniques from semigroup theory.

First, since $u(t) \in H^1_0(0,L)$, from the equation \eqref{160} we deduce that
\[
\aligned
\Imag \langle i\del_t u, \overline u \rangle_{H^{-1} \times H^1_0} + 
\Imag \langle \del_{xx} u , \overline u \rangle_{H^{-1} \times H^1_0} =0,
\endaligned
\]
and so,
\[
\aligned
\frac{d}{dt} \int_0^L |u|^2(t) \,dx = 0,
\endaligned
\]
which gives \eqref{250}.

For the proof of \eqref{260}, we follow the ideas of Kato \cite[Lemma 3.1]{18}.
We point out that $-\del^3_x$ is a skew-adjoint operator on $L^2(0,L)$ with domain 
$H^3(0,L) \cap H^2_0(0,L)$ and it generates a group of isometries $U_K(t)$ on $L^2(0,L)$. Also,
$i\Delta$ is a skew-adjoint operator on $L^2(0,L)$ with domain 
$H^2(0,L) \cap H^1_0(0,L)$ and so it generates a group of isometries $U_S(t)$ on $L^2(0,L).$ 
Now, we write the equations \eqref{155} in integral form,
\be
\label{270}
u(t) = U_S (t-s) u(s) + \int_s^t U_S (t - r) m(r) \,dr,
\ee
\be
\label{280}
v(t) = U_K (t-s) v(s) + \int_s^t U_K (t - r) n(r) \,dr,
\ee
with
\[
\aligned
&m(r) = -i\alpha g^M (v(r)) u(r) - i\beta |u(r)|^2 u(r),
\\
& n(r) = -\del_x f^M(v(r)) + \gamma \del_x \big( (g^M)' (v(r)) |u(r)|^2 \big).
\endaligned
\]
The formulas \eqref{270},\eqref{280} are verified in $L^2(0,L)$ and, since $\del_x m
\in L^\infty_\loc (I; L^2(0,L))$, we have also
\be
\label{290}
\del_x u(t) = U_S (t-s)\del_x u(s) + \int_s^t U_S (t - r)\del_x m(r) \,dr
\ee
in $L^2(0,L)$.

Using the isometric property of $U_K(t)$ and $U_S(t)$, we easily deduce from \eqref{270}--\eqref{290}
\[
\big( v(t) , v(t) \big) = \big( v(s) , v(s) \big) + 2 \int_s^t \big( v(r) , n(r) \big) \,dr,
\]
\[
\Imag \big( u(t), \del_x u(t) \big) = \Imag \big (u(s) , \del_x u(s) \big) + 2 \int_s^t \Imag \big( m(r),
\del_x u(r) \big) \,dr.
\]
Since $u(r), v(r)\in H^1_0(0,L)$, we obtain from the expressions of $m(r)$ and $n(r)$
\[
\aligned
\Imag \big(m(r), \del_x u(r) \big) &= - \Real \alpha \int_0^L g^M(v(r)) u(r) \del_x \overline  u(r) \,dx
\\
&\qquad - \beta \Real \int_0^L | u(r)|^2 u(r) \del_x \overline  u(r) \,dx 
\\
& = \frac\alpha2 \int_0^L (g^M)'(v(r)) \del_x v(r) |u(r)|^2 \,dx.
\endaligned
\]
It follows that
\[
\aligned
2\gamma \Imag \big( m(r) , \del_x u(r) \big) = - \alpha \big( v(r), n(r) \big),
\endaligned
\]
which implies the conclusion \eqref{260}. This completes the proof of Proposition~\ref{245}.
\endproof

The following result establishes the convergence of the approximations \eqref{185} to a global solution
of the truncated problem \eqref{155}, and the crucial energy estimate \eqref{330}.

\begin{proposition}
\label{305}
Let $\alpha, \beta, \gamma$ be such that $\alpha \gamma>0$. For each $M>0$ let $(u^{h,M}, v^{h,M})
\in (C([0,T]; X_J))^2$ be the solution of \eqref{185} with initial data $(u^h_0, v^h_0)$ such that
$\PP_1^h u^h_0 \to u_0$, $ \PP_1^h v^h_0 \to v_0$ in $H^1(0,L)$ as $h\to 0$. Then, up to a 
subsequence, 
\be
\label{300}
\PP_1^h u^{h,M} \ws u^M, \quad \PP_1^h v^{h,M} \ws v^M \qquad \text{ in }\quad
L^\infty\big( [0,T] ; H^1(0,L) \big),
\ee
\be
\label{310}
\PP_1^h u^{h,M} \to u^M, \quad \PP_1^h v^{h,M} \to v^M \qquad \text{ in }\quad
L^\infty\big( [0,T] ; L^2(0,L) \big),
\ee
with 
\be
\label{320}
\aligned
(u^M, v^M) \in {}&\big( L^\infty ([0,T]; H^1_0(0,L)) \big)^2 
\\
&\cap \big( C( [0,T] ; H^1_0(0,L)) \times 
C( [0,T] ; L^2(0,L)) \big), \quad T>0,
\endaligned
\ee
a global weak solution of the truncated
system \eqref{155}. Moreover, the following energy estimate is valid,
\be
\label{330}
\aligned
\Ecal^M (t)  :=  \int_0^L \big\{ \gamma |\del_x u^M(t) |^2 {}&+ \frac\alpha2 |\del_x v^M(t)|^2
+ \alpha\gamma g^M(v^M(t)) |u^M(t)|^2 
\\
&\quad- \alpha F^M (v^M(t)) + \frac{\beta\gamma}2 |u^M(t)|^4 \big\} \,dx
\le \Ecal^M(0),
\endaligned
\ee
for all $t \in [0,T].$
\end{proposition}

We postpone the proof of Proposition \ref{305} until Section~\ref{335}, and proceed to prove 
Theorem~\ref{225}. The goal is to prove an $L^\infty$ bound for $u^M,v^M$ independent of the
truncation parameter $M$, using the energy inequality \eqref{330}. Once this is achieved, it is clear from the definition of $f^M$ and $g^M$ 
that taking $M$ large enough yields a solution of the original problem,~\eqref{15}.

Let us define $\Mcal_0 = \Mcal(0)$, $\Qcal_0 = \Qcal(0)$ (see 
\eqref{250},\eqref{260}), and set
\be
\label{331}
\aligned
\Ecal_0 = |\gamma| \| \del_x u_0 \|^2_2 + \frac{|\alpha|}2 \| \del_x v_0\|_2^2 + |\alpha \gamma| 
\| v_0 \|_2 \| u_0 \|_4^2 + \frac{|\alpha|}3 \| v_0 \|_3^3 + \frac{|\beta\gamma|}2 \|u_0 \|_4^4.
\endaligned
\ee
Observe that $|\Ecal^M(0)| \le \Ecal_0$ for all $M>0$ (see \eqref{330}). 

For simplicity, in what follows we omit the superscript $M$ from the solutions $(u^M, v^M)$ of 
the system \eqref{155} obtained in Proposition~\ref{305}. 

First, note that the energy inequality \eqref{330} gives
\be
\label{331.10}
\aligned
|\gamma| \| \del_x u(t) \|_2^2 +{}& \frac{|\alpha|}2\| \del_x v(t) \|_2^2 
\\
&\le  \big( \Ecal_0 + |\alpha\gamma| \| v(t) \|_2 \| u(t) \|_4^2 + \frac{|\alpha|}3 
\| v(t) \|_3^3 + \frac{|\beta\gamma|}2 \| u(t) \|_4^4 \big).
\endaligned
\ee
Next, from \eqref{260} we have
\be
\label{332}
\aligned
\| v(t) \|_2^2 \le \frac1{|\alpha|} \big( | \Qcal_0| + 2 | \gamma| \|u_0\|_2\| \del_x u(t) \|_2 \big).
\endaligned
\ee
Let now $m = \min\{|\gamma| , |\alpha|/2\}$. Using again  Gagliardo--Nirenberg and Young 
inequalities, we deduce from \eqref{331.10}--\eqref{332} (as in \cite{13})
\[
\aligned
\| \del_x u(t) \|_2^2 + \| \del_x v(t) \|_2^2 &\le \frac1m \big( |\gamma|
\| \del_x u(t) \|_2^2 + \frac{|\alpha|}2 \| \del_x v(t) \|_2^2 \big)
\\
&\le \frac1m \big( \Ecal_0 + |\alpha\gamma| \| v(t) \|_2 \| u(t) \|_4^2 + \frac{|\alpha|}3 
\| v(t) \|_3^3 + \frac{|\beta\gamma|}2 \| u(t) \|_4^4 \big)
\\
& \le  C\big( \Ecal_0 + \| v(t) \|_2^2 + \| v(t) \|_3^3 + \| u(t)\|_4^4 \big) 
\\
&\le C \big( \Ecal_0 + |\Qcal_0| + |\Qcal_0|^{5/3} + \Mcal_0 + \Mcal_0^3 + \Mcal_0^5 \big)
\\
&\qquad+ \frac12( \| \del_x u(t) \|_2^2 + \| \del_x v(t) \|_2^2 ),
\endaligned
\]
with $C = C(\alpha,\beta, \gamma) $ only depending on the parameters $\alpha,\beta,\gamma$.
Therefore
\[
\aligned
 \| \del_x u(t) \|_2^2 + \| \del_x v(t) \|_2^2 
\le 2C \big( \Ecal_0 + |\Qcal_0| + |\Qcal_0|^{5/3} + \Mcal_0 + \Mcal_0^3 + \Mcal_0^5 \big)
 :=  K_0.
\endaligned
\]
Note that $K_0$ is independent of $M$.
Finally, since from \eqref{332},
\[
\aligned
\| v(t) \|_\infty^2 &\le 2 \| v(t) \|_2  \| \del_x v(t) \|_2 \le \| v(t) \|_2^2 + \| \del_x v(t) \|_2^2
\\
& \le \frac1{|\alpha|} |\Qcal_0| + \Big| \frac\gamma\alpha \Big| \Big(  \| u_0 \|_2^2 + \| \del_x u(t) \|_2^2 \Big)
+\| \del_x v(t) \|_2^2,
\endaligned
\]
we obtain
\be
\label{334}
\| v(t) \|_\infty \le \big( \frac1{|\alpha|} | \Qcal_0| + \Big| \frac\gamma\alpha \Big| 
\Mcal_0 + \big(1+ |\gamma/\alpha| \big) K_0 \big)^{1/2}  :=  \overline K,
\ee
with $\overline K$ independent of $M$ but depending on the initial data.
Therefore, if the truncation parameter $M$ in \eqref{155} satisfies $M > \overline  K$, we conclude by
\eqref{334} and the definition of  $f^M$, $g^M$ that $(u, v)  :=  (u^M, v^M)$ is actually a 
solution of the Sch\-r\"o\-din\-ger--KdV system \eqref{15}. This completes the proof of 
Theorem~\ref{225}.

\section{Proof of Proposition \ref{305}}
\label{335}
First of all, we need the following lemma concerning the global existence of the solution of
the discrete problem \eqref{185}.
Due to the fact that the problem \eqref{185} is truncated, we are also able to obtain the essential
$H^1$ estimate \eqref{340}. For simplicity, we will omit the superscript $M$.
\begin{lemma}
\label{336}
Let $\alpha,\beta,\gamma \in \RR$ be such that $\alpha\gamma >0$. Fix $J\in \NN$, $L>0$, 
$h= L/(J+1)$, and $(u^h_0,v^h_0)
\in X_{J,\CC} \times X_{J, \RR}$. Then, for each $T>0$, there exists a unique solution 
\[
\aligned
(\uh, \vh) \in C( [0,T]; X_{J,\CC} ) \times C( [0,T]; X_{J,\RR} )  
\endaligned
\]
of the problem \eqref{185}. Moreover, the following estimate is valid,
\be
\label{340}
\aligned
\| D_+ \uh (t) \|_2^2 {}&+ \| D_+ \vh (t) \|_2^2 + \| \uh (t) \|_2^2 + \| \vh (t) \|_2^2 
\\
& \le C \big( \| u^h_0\|_2, \| v^h_0\|_2, \|D_+ u^h_0\|_2, \|D_+ v^h_0\|_2, T, M),
\endaligned
\ee
for all $t\in [0,T]$.
\end{lemma}
\proof
Let $S_h (t) = e^{i\Delta^ht},$ $G_h (t) = e^{-D^3 t}$ be the unitary groups generated by the 
discrete operators $i\Delta^h$ and $-D^3$ in the $X_J$ space. The problem \eqref{185} can be written
in the usual semigroup framework, as an integral equation in the $ X_{J,\CC} \times X_{J, \RR}$ space:
\be
\label{350}
\aligned
 \uh(t) & = S_h (t) u^h_0  
+ \int_0^t  S_h( t - s) J_S(\uh(s) , \vh(s)) \,ds 
=:  \Phi_1(\uh, \vh),
\\
 \vh(t) & =  G_h (t) v^h_0
+ \int_0^t  G_h( t - s) J_K(\uh(s) , \vh(s)) \,ds 
=:  \Phi_2(\uh, \vh),
\endaligned
\ee
where
\[
J_S (\uh, \vh) = -i\alpha g(\vh) \uh - i \beta |\uh|^2 \uh,
\]
\[
J_K (\uh, \vh) = - D_0 f(\vh) + \gamma D_0 (g'(\vh) |\uh|^2).
\]
For $R> \| u^h_0 \|_2 + \| v^h_0 \|_2 $ and $T>0$ we consider the product space 
$B^T_{R,\CC} \times B^T_{R,\RR},$ with
\[
B^T_{R,\CC} = \big\{ w \in C([0,T] ; X_{J,\CC}) : \|w\|_{L^\infty([0,T]; X_J)} \le R \big \}
\]
and similarly for $B^T_{R,\RR}$.

By \eqref{100}, it is now a simple matter to prove that there exists $T>0$ such that the map
\[
(\uh, \vh) \in B^T_{R,\CC} \times B^T_{R,\RR} \mapsto \Phi( \uh, \vh) := (\Phi_1( \uh, \vh),
 \Phi_2( \uh, \vh) )
\]
is a strict contraction on the complete metric space $B^T_{R,\CC} \times B^T_{R,\RR}$. Thus, by the
Banach fixed-point theorem, we obtain a unique local in time solution $(\uh, \vh)$ of the problem
\eqref{185} in the space $ C( [0,T]; X_{J,\CC} ) \times C( [0,T]; X_{J,\RR} ).$

To obtain a global solution, we must estimate the $l^2$-norm of $\uh(t)$ and $\vh(t)$, for each 
fixed $h$. 
From the equation \eqref{190}, the conservation of the $l^2$-norm of $\uh(t)$ follows easily
by taking the scalar product with $\uh$ and summation by parts. Applying the same 
procedure to equation \eqref{200}, we find
\be
\label{360}
\aligned
\frac12\del_t \| \vh(t) \|_2^2 = \sum_{j=1}^{J} h f(\vj) D_0 \vj - \sum_{j=1}^{J} hg'(\vj) |\uj|^2 D_0 \vj.
\endaligned
\ee
From the definition of $f,g$, the conservation of the $l^2$-norm of $\uh$ and the fact that, 
for $h$ fixed, $\| \cdot \|_\infty \le C(h) \| \cdot \|_2$, we derive that
\[
\| \vh(t) \|_2^2 \le \| v^h_0 \|_2^2 + C(h) \int_0^t \| \vh(s) \|_2^2 \,ds.
\]
The conclusion now follows from a Gronwall argument.

It remains to prove the inequality \eqref{340}. In addition to the conservation of the $l^2$-norm of
$\uh$, we have the conservation of the discrete energy:
\be
\label{370}
\aligned
E^h(t) & :=  \gamma \|D_+ \uh(t) \|_2^2 + \frac\alpha2 \|D_+ \vh(t) \|_2^2 +
\frac{\beta\gamma}2 \|\uh(t) \|_4^4 
\\
&\quad+ \alpha\gamma \sum_{j=1}^{J} h g(\vj) |\uj|^2 - \alpha \sum_{j=1}^{J} h F(\vj) = E^h(0).
\endaligned
\ee
To prove this identity, we proceed in the same way as in the continuous case \cite{11}: Take the 
scalar product in $X_J$ of the equation \eqref{190} with $\del_t \overline u^h$, take the 
real part, and use the equation \eqref{200} and the skew-adjoint properties of the operators 
$D_0$ and $D^3$.

Now we return to \eqref{360} and observe that from $f(\vj) = f(\vj) - f(0) = f'(\theta_j) \vj,$ $|g'|\le1$
and  \eqref{90} we find
\[
\aligned
\del_t \| \vh \|_2^2 &\le C(M) \big( \| D_+ \vh \|_2^2 + \| \vh \|_2^2 \big) + C\| \uh\|_\infty 
\| D_+ \vh\|_2
\\
& \le C(M) \big( \| D_+ \vh \|_2^2 + \| D_+ \uh\|_2^2 + \| \vh \|_2^2 \big). 
\endaligned
\]
Integrating on $(0,t)$ and using Gronwall's lemma, we obtain
\be
\label{380}
\| \vh (t)\|_2 \le a_1(t,M) + a_2(t,M) \int_0^t  \| D_+ \vh(s) \|_2^2 + \| D_+ \uh(s)\|_2^2 \,ds
\ee
for some continuous functions $a_1,a_2$. On the other hand, from the conservation of the energy
\eqref{370} and using \eqref{80} we get
\be
\label{390}
\aligned
\| D_+ \uh(s)\|_2^2 + \| D_+ \vh(s)\|_2^2 &\le C(u^h_0, v^h_0) + C_1 \| D_+ \uh(s)\|_2
\\
&\quad+ C_2 \sum_{j=1}^{J} h g(\vj)|\uj|^2 + C_3 \sum_{j=1}^{J} h F(\vj).
\endaligned
\ee
But now, the definition of $f^M$ allows us (roughly) to bound $F^M(v)$ by $C(M) v^2$. This is essential 
in view of the desired $H^1$ estimate \eqref{340}, since these terms would otherwise be cubic.
We have
\[
\aligned
\sum_{j=1}^{J} h |F(\vj)| &= \sum_{j=1}^{J} h |F(\vj) - F(0)| = \sum_{j=1}^{J} h |f(\theta_j) \vj|
\\
&\le \sum_{j=1}^{J} h f(\theta_j)^2 + \| \vh\|^2_2,
\endaligned
\]
for some $\theta_j$ between 0 and $\vj$. Now,
\[
\aligned
\sum_{j=1}^{J} h f(\theta_j)^2 = \sum_{|\vj| \le M^2+1} hf(\theta_j)^2 + 
\sum_{|\vj| > M^2+1} hf(\theta_j)^2. 
\endaligned
\]
Recall the definition of the truncated functions in \eqref{150}. 
For the first sum, we have $ f(\theta_j)^2 \le \vj^4 \le (M^2+1)^2 \vj^2$, and for the second
sum we have $ f(\theta_j)^2 \le \vj^2$. Thus we obtain
\[
\sum_{j=1}^{J} h |F(\vj)| \le C(M) \|\vh \|_2^2.
\]
Similarly, since the $l^2$-norm of $\uh$ is conserved, we find
\[
\sum_{j=1}^{J} h | g(\vj)| |\uj|^2 \le \| g\|_\infty \| \uh\|_2^2 \le C(M) \|u^h_0\|_2^2.
\]
These estimates together with \eqref{380} and \eqref{390} give us
\[
\aligned
\| D_+ \uh(t)\|_2^2 + \| D_+ \vh(t)\|_2^2 &\le C(u^h_0, v^h_0, M) + C(M) \| \vh (t) \|_2^2
\\
&\hspace{-15mm}\le C(u^h_0, v^h_0, M) + C(t,M) \int_0^t \| D_+ \uh(s)\|_2^2 + \| D_+ \vh(s)\|_2^2 \,ds,
\endaligned
\]
where $C(u^h_0, v^h_0, M) = 
C(\| u^h_0 \|_2, \| D_+u^h_0 \|_2, \| v^h_0 \|_2, \|D_+ v^h_0 \|_2, M)$. A Gronwall argument,
\eqref{380}, and the conservation of $\|\uh\|_2$ give the conclusion \eqref{340}. This completes 
the proof of Lemma~\ref{336}.
\endproof

\proof[Proof of Proposition~\ref{305}]
We will use the interpolators $ \PP_1^h, \PP_0^h$ defined 
in \eqref{50}. Since $\del_x \PP_1^h = \PP_0^h D_+$, it follows from the hypothesis 
$ \PP_1^h u^h_0 \to u_0,\ \PP_1^h v^h_0 \to v_0$ in $H^1_0(0,L)$, \eqref{70} and \eqref{340} that
\[
\| \PP_1^h \uh (t)\|_{H^1(0,L)} \le C, \quad \| \PP_1^h \vh (t)\|_{H^1(0,L)} \le C,\quad 
t\in [0,T]
\]
with $C = C( \|u_0\|_{H^1}, \|v_0\|_{H^1}, T, M).$ Thus, using the compactness of the embedding
of $H^1(0,L)$ into $L^2(0,L)$, we obtain, as $h\to 0$ (for a subsequence),
\be
\label{400}
\aligned
&\PP_1^h \uh \ws u, \quad \PP_1^h \vh \ws v \qquad \text{ in }\quad
L^\infty\big( [0,T] ; H^1(0,L) \big),
\\
&\PP_1^h \uh \to u, \quad \PP_1^h \vh \to v \qquad \text{ in }\quad
L^\infty\big( [0,T] ; L^2(0,L) \big),
\endaligned
\ee
for some $u,v \in H^1_0(0,L)$.
Also, we have from \eqref{60},
\be
\label{410}
\PP_0^h \uh \to u, \quad \PP_0^h \vh \to v \text{ in } L^\infty\big( [0,T] ; L^2(0,L) \big).
\ee
To prove that $u,v$ are solutions to the system \eqref{155}, we apply the piecewise constant 
interpolator $ \PP_0^h$ to the equations \eqref{190},\eqref{200}:
\be
\label{420}
i\del_t \PP_0^h \uh + \PP_0^h \Delta^h \uh = \beta \PP_0^h(|\uh|^2 \uh) 
+\alpha \PP_0^h(g^M(\vh) \uh)
\ee
\be
\label{430}
\del_t \PP_0^h \vh +  \PP_0^h D^3\vh +  \PP_0^h D_0 f^M ({\vh}) 
= \gamma \PP_0^h D_0( (g^M)'(\vh) | \uh|^2).
\ee
From \eqref{60},\eqref{340} we have
\be
\label{440}
\PP_1^h f(\vh) - \PP_0^h f(\vh) \to 0 \text{ in } L^\infty\big( [0,T] ; L^2(0,L) \big)
\ee
and, since the piecewise constant interpolator commutes with nonliearities, it follows from \eqref{410}
that
\be
\label{450}
\PP_0^h f(\vh) = f (\PP_0^h \vh)\to f(v)\text{ in } L^\infty\big( [0,T] ; L^2(0,L) \big).
\ee
On the other hand, using \eqref{340},
\[
\| \del_x \PP_1^h f(\vh) \|_2 = \|  \PP_0^h D_+f(\vh) \|_2 = \| D_+ f(\vh) \|_2 \le C(M),
\]
and, from \eqref{440},\eqref{450} we deduce
\[
\PP_1^h f(\vh) \ws  f (v) \text{ in } L^\infty\big( [0,T] ; H^1(0,L) \big)
\]
and so,
\[
\del_x \PP_1^h f(\vh) \ws  \del_x f (v) \text{ in } L^\infty\big( [0,T] ; L^2(0,L) \big).
\]
Similarly, and using also \eqref{410}, we find
\[
\aligned
&\PP_0^h (| \uh|^2 \uh) = | \PP_0^h \uh|^2 \PP_0^h\uh  \to  |u|^2 u 
\text{ in } L^\infty\big( [0,T] ; L^2(0,L) \big),
\\
&\PP_0^h (g(\vh) \uh) = g( \PP_0^h \vh) \PP_0^h\uh  \to  g(v) u 
\text{ in } L^\infty\big( [0,T] ; L^2(0,L) \big),
\\
&\PP_0^h D_0 (g'(\vh) |\uh|^2) = \del_x \PP_1^h (g'(\vh) |\uh|^2) \ws \del_x (g'(v) |u|^2)
\text{ in } L^\infty\big( [0,T] ; L^2(0,L) \big),
\endaligned
\]
which allows us to pass to the limit on the corresponding terms in the weak formulation
of the equations \eqref{15}.

It remains to analyze the terms $ \PP_0^h \Delta^h \uh$  and $ \PP_0^h D^3 \vh$. Let 
$\phi \in \Dcal(0,L)$ be a test function. We have
\[
\aligned
\big \langle \PP_0^h D^3 \vh , \phi \big \rangle & = \sum_{j=2}^{J-1} \int_{\xj}^{\xjj} \PP_0^h
D^3 \vh \phi \,dx
 = \sum_{j=2}^{J-1} D_0D_-D_+\vj \int_{\xj}^{\xjj} \phi(x) \,dx
\\
& = \sum_{j=2}^{J-1} D_+ \vj \int_{\xj}^{\xjj} \frac1{2h^2} \big( \phi(x-2h) - \phi(x-h) + \phi(x) 
-\phi(x+h) \big) \,dx
\endaligned
\]
and so, by Taylor expansion of $\phi$,
\[
\aligned
\big | \big \langle \PP_0^h D^3 \vh , \phi \big\rangle \big| &\le C \sum_{j=2}^{J-1}
h | D_+ \vj| \| \phi''\|_\infty \le C \Big( \sum_{j=2}^{J-1}
h | D_+ \vj|^2 \Big)^{1/2}  \| \phi\|_{H^3(0,L)}.
\endaligned
\]
Hence, from \eqref{340} we obtain
\be
\label{460}
\| \PP_0^h D^3 \vh \|_{ L^\infty( [0,T] ; H^{-3}(0,L))} \le C.
\ee
If we now take a test function $\varphi \in \Dcal ((0,T)\times(0,L)),$ we may compute in the sense of
distributions
\[
\aligned
\big \langle \PP_0^h D^3 \vh , \varphi \big\rangle & = \int_{0}^T \sum_{j=1}^{J} D^3 \vj 
\int_{\xj}^{\xjj} \varphi (x,t) \,dx \,dt
\\
&= - \int_{0}^T \sum_{j=1}^{J} \vj \int_{\xj}^{\xjj} ( \del^3_x \varphi + \OO (h) ) \,dx \,dt
\\
&= - \big\langle \PP_0^h \vh, \del^3_x \varphi \rangle + \OO(h) \to
-\langle v, \del^3_x \varphi \rangle = \langle  \del^3_x v,  \varphi \rangle
\endaligned
\]
as $h\to 0$. Hence, we deduce from \eqref{460} that
\[
\PP_0^h D^3 \vh \ws  \del^3_x v  \text{ in } L^\infty\big( [0,T] ; H^{-3}(0,L) \big).
\]
In a similar way we prove that
\[
\PP_0^h \Delta^h \uh \ws  \Delta u  \text{ in } L^\infty\big( [0,T] ; H^{-2}(0,L) \big)
\]
and using the equations,
\[
\aligned
&i \del_t\PP_0^h \uh \ws  i \del_t u  \text{ in } L^\infty\big( [0,T] ; H^{-2}(0,L) \big)
\\
&\del_t \PP_0^h \vh \ws  \del_t v  \text{ in } L^\infty\big( [0,T] ; H^{-3}(0,L) \big).
\endaligned
\]
Therefore, taking the limit $h\to0$ in the weak formulation of the
equations \eqref{420},\eqref{430} we obtain a weak solution
\[
\aligned
(u, v) \in {}&\big( L^\infty ([0,T]; H^1_0(0,L)) \big)^2 
\\
&\cap \big( C( [0,T] ; H^{-2}(0,L)) \times 
C( [0,T] ; H^{-3}(0,L)) \big), \quad T>0,
\endaligned
\]
of the problem \eqref{155}. To prove \eqref{320}, recall that this solution satisfies the integral
system \eqref{270},\eqref{280}. Since 
\[
\| m(u,v) \|_{H^1} \le C ( \|u_0\|_{H^1}, \|v_0\|_{H^1} ),
\]
\[
\| n(u,v) \|_{H^1} \le C ( \|u_0\|_{H^1}, \|v_0\|_{H^1} ),
\]
we deduce from \eqref{460} that
\[
(u,v) \in C([0,T]; H^1) \times C([0,T];L^2).
\]
It remains to prove the energy inequality \eqref{330}. Let us write the discrete energy \eqref{370} 
in the form
\be
\label{470}
\aligned
E^h(t) & :=  \gamma \|\del_x \PP_1^h \uh(t) \|_2^2 + \frac\alpha2 \|\del_x \PP_1^h\vh(t) \|_2^2 +
\frac{\beta\gamma}2 \| \PP_0^h \uh(t) \|_4^4 
\\
&\quad+ \alpha\gamma \int_0^L  g( \PP_0^h\vh) | \PP_0^h\uh|^2 \,dx - \alpha \int_0^L  F( \PP_0^h\vh) \,dx = E^h(0).
\endaligned
\ee
Now we recall the weak and strong convergences \eqref{400},\eqref{410}. From the last term on the
left-hand side of \eqref{470}, and since $|F'(\xi)| \le C|\xi|^2$, we find
\[
\aligned
\int_0^L \big | F( \PP_0^h \vh(t) ) - F(v(t)) \big| \,dx &\le C \| \PP_0^h \vh(t) - v(t) \|_2
\| (\PP_0^h \vh(t))^2 + v^2(t) \|_2
\\
& \le C  \| \PP_0^h \vh(t) - v(t) \|_2 \to 0 \quad(h\to0)
\endaligned
\]
and so 
\[
\int_0^L F( \PP_0^h \vh(t)) \,dx \to \int_0^L F(v(t)) \,dx.
\]
Note that here it is essential that the spatial domain $(0,L)$ be bounded. Indeed, a version 
of the energy inequality \eqref{330} on the whole line cannot be obtained using the 
available convergences \eqref{400},\eqref{410}, which are local in space.

Finally, from the strong convergence in $H^1(0,L)$ of the initial data $(u^h_0, v^h_0)$ 
and using the lower semi-continuity of the $H^1$ norm, we easily obtain from \eqref{470},
in the limit $h\to0$, the conclusion \eqref{330}: $\Ecal(t) \le \Ecal(0),$ $t\in[0,T]$.
This completes the proof of Proposition~\ref{305}. \endproof


\section{Numerical experiments}
In this section, we present some numerical computations using a fully discrete version of the
method \eqref{115}. We emphasize that these simulations are for the sake of illustration of our convergence result only.
In particular, it would be interesting to perform more extensive numerical tests, such as determining the order of convergence, 
or employing more sophisticated time discretizations,
which we 
do not perform here.

Given some time step $\tau>0$, a spatial mesh size $h$, and initial data
$(u_{0j}, v_{0j})_{j=0,\dots,J+1}$, we consider for $n\ge 0$ 
the following algorithm. Set $u^{n+1/2}_j = \frac12(\unnj + \unj)$ and solve
\be
\label{scheme}
\aligned
& i \frac1\tau( \unnj - \unj) +  \Delta^h u^{n+1/2}_j
= \big| u^{n+1/2}_j \big|^2 u^{n+1/2}_j + \vnj u^{n+1/2}_j,
\\
&\frac1\tau(\vnnj - \vnj) +  D^3 \vnnj  
+ D_0 (\vnnj)^{2} = D_0( | \unj|^2),
\endaligned
\ee
for $j=0,\dots, J+1$. This corresponds to a semi-implicit Crank--Nickolson scheme for the Sch\-r\"o\-din\-ger equation and a 
fully implicit Euler scheme for the KdV equation. Because of the nonlinear terms, we perform a Newton 
iteration at each time step with a tolerance of $10^{-6}$. At each Newton iteration we solve 
independently a tridiagonal system for the first equation of \eqref{scheme} by a standard direct 
method, and the pentadiagonal system issuing from the second equation by an $LU$ 
decomposition method.

\subsection{Comparison with exact solutions}
We now test our scheme and illustrate our convergence result. We will simulate the following system,
\be
\label{500}
\left\{
\aligned
&i\del_t u + \del_{xx} u = \alpha v u - |u|^2 u 
\\
& \del_t v +  \del_{x}^3 v + v \del_x v = \frac\alpha2 \del_x| u|^2,
\endaligned\right.
\ee
which is the same as \eqref{15} for a special choice of the parameters, except for the quasilinear 
term in the KdV equation which is formally equivalent to $\frac12 \del_x(v^2)$. In \cite{DFO} 
we can find the following exact 
traveling wave solutions to \eqref{500},
\[
\aligned
(u,v) = ( e^{i\omega t}e^{ixc/2} \phi(x-ct), \psi (x-ct)),
\endaligned
\]
with 
\[
\aligned
\phi(y) = \frac{\sqrt{2c^* (1+6\alpha)}}{\cosh (\sqrt{c^*}y)}, \qquad \psi (y) 
= \frac{12 c^*}{\cosh^2(\sqrt{c^*}y)}.
\endaligned
\]
Here, $\alpha \in [-1/6, 0]$ and $\omega \in\RR$ are given, and
$2c = 1 + \sqrt{1+\frac\alpha3(1+6\alpha)},$ $c^* = c^2/4 + \omega^2$.
We chose $\alpha = -1/12$ and $\omega = 0$ for the simulations below. This gives a traveling wave
speed $c=0.996516$.

\begin{figure}
\includegraphics[width=\linewidth,keepaspectratio=false]{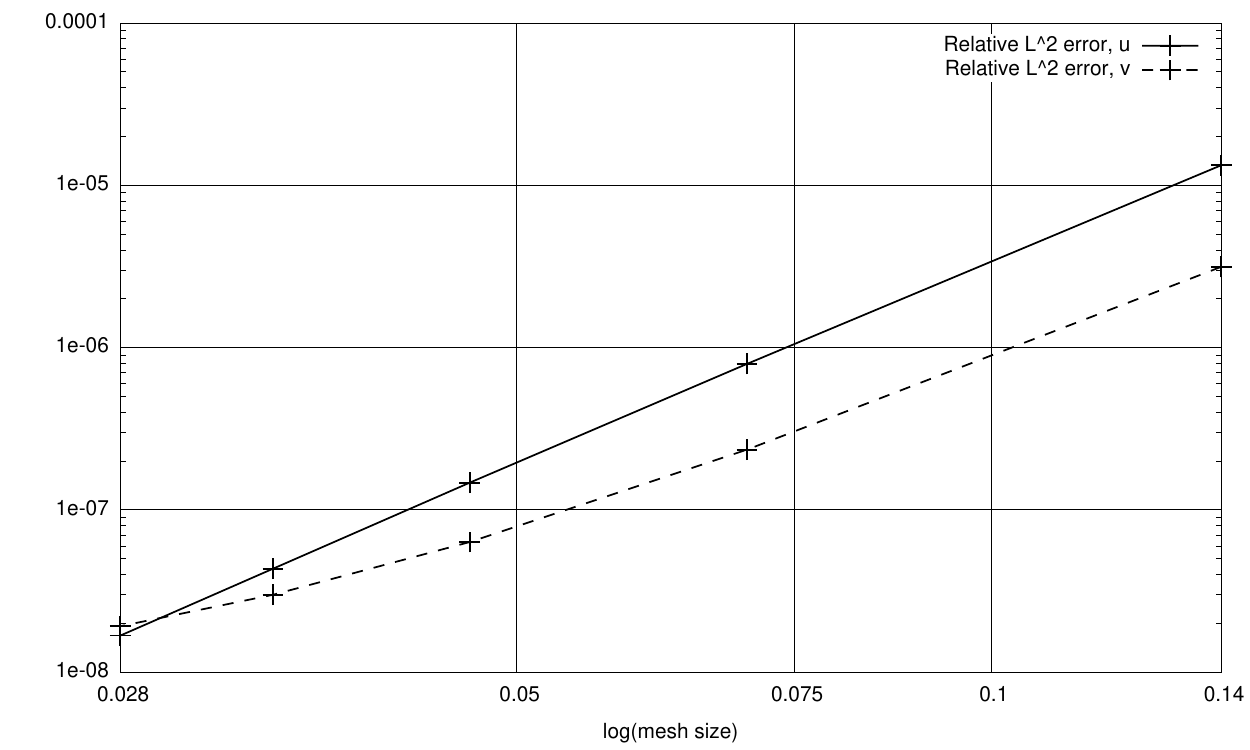}
\caption{Relative $L^2$ error, $T=5$, $\tau = 0.0001$, as a function of mesh size. 500 to 2500 spatial points.}
\label{550}
\end{figure}

\begin{figure}
\includegraphics[width=\linewidth,keepaspectratio=false]{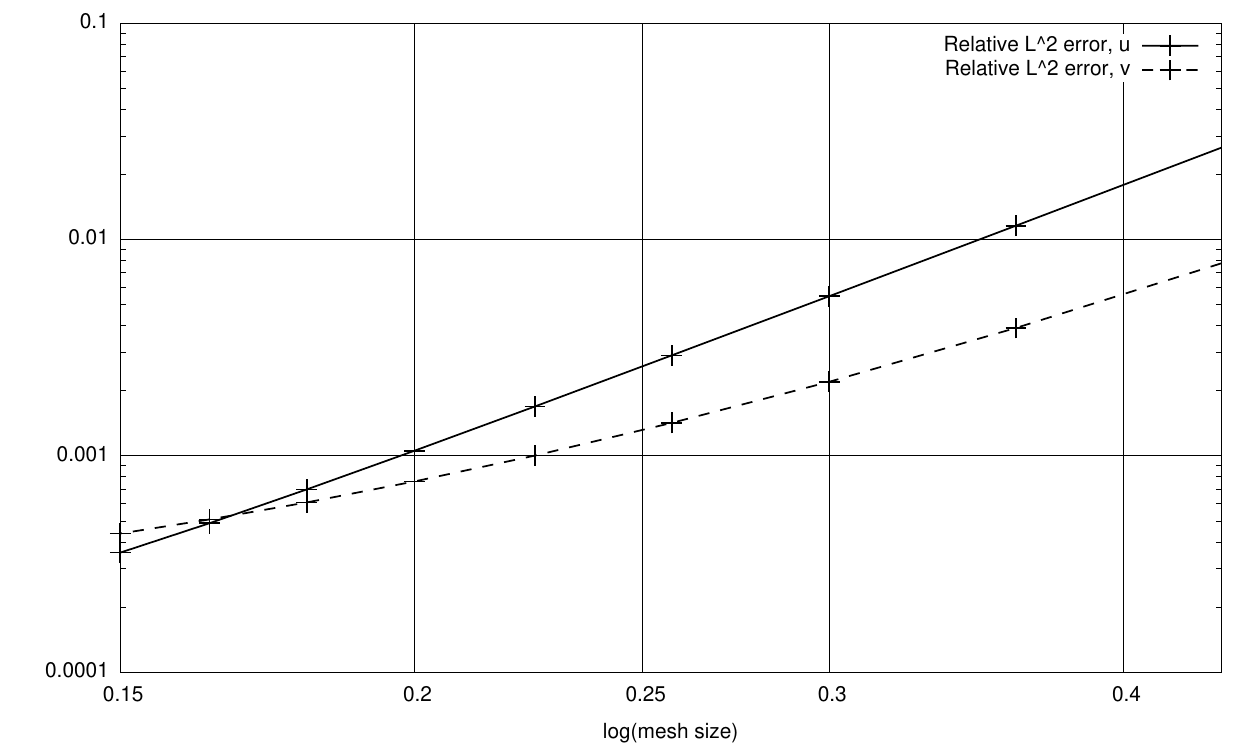}
\caption{Relative $L^2$ error, $T=30$, $\tau = 0.0005$, as a function of mesh size. 200 to 600 spatial points.} 
\label{560}
\end{figure}

In Figure \ref{550} we present (in logarithmic scale) the relative $L^2$ error computed at $T = 5$ for a time step 
$\tau = 0.0001$ as a function of the mesh size. The computational domain is $[-20,50]$ and the number of spatial points ranges from 500 to 2500.
In Figure~\ref{560} we present the relative $L^2$ error for a similar computation with $T=30$ 
and a time step $\tau =0.0005$. The number of spatial points ranges from 200 to 600.

In can be seen that the exact solution is approximated rather well by our simple numerical scheme, 
especially bearing in mind that the simulations (which, again, serve only to illustrate our results) were performed in a few minutes
on a laptop running at 2.4 GHz.

\section*{Acknowledgements}
The authors were supported by the Portuguese Foundation for Science and Technology (FCT), PEst OE/MAT/UI0209/2011, and the FCT grant PTDC/MAT/110613/2009.
PA was also supported by FCT through a Ci\^encia 2008 fellowship.


%



\begin{thebibliography}{}
%
%


\bibitem{12}{D. Bekiranov, T. Ogawa and G. Ponce,} Interaction equations for short and long dispersive waves, J. Funct. Anal. 158 (1998), no. 2, 357--388

\bibitem{Benney}{D.J. Benney,}
{A general theory for interactions between short and long waves,}
Stud. Appl. Math. 56 (1977) 81--94

\bibitem{14} J. Bourgain,
Fourier transform restriction phenomena for certain lattice subsets and applications to nonlinear evolution equations. I. Schr\"odinger equations. Geom. Funct. Anal. 3 (1993), no. 2, 107--156.

\bibitem{13} A.J. Corcho and F. Linares,
Well-posedness for the Schršdinger-Korteweg-de Vries system. Trans. Amer. Math. Soc. 359 (2007), no. 9, 4089--4106.

\bibitem{17}  {J.-P. Dias, M. Figueira and F. Oliveira,}
Well-posedness and existence of bound states for a coupled Sch\-r\"o\-din\-ger--gKdV system,
Nonlinear Anal. 73 (2010), no. 8, 2686--2698.

\bibitem{DFO}  {J.-P. Dias, M. Figueira and F. Oliveira,}
existence of bound states for the coupled Sch\-r\"o\-din\-ger--gKdV system with cubic nonlinearity,
C. R. Math. Acad. Sci. Paris 384 (2010), no. 19-20, 1079--1082.

\bibitem{18} T. Kato,
On the Cauchy problem for the (generalized) Korteweg-de Vries equation. Studies in applied mathematics, 93--128, 
Adv. Math. Suppl. Stud., 8, Academic Press, New York, 1983.

\bibitem{8}
T. Kawahara, N. Sugimoto and T. Kakutani,
Nonlinear Interaction between Short and Long Capillary-Gravity Waves,
J. Phys. Soc. Japan 39 (1975), 1379--1386.

\bibitem{15} C.E. Kenig, G. Ponce and L. Vega,
The Cauchy problem for the Korteweg-de Vries equation in Sobolev spaces of negative indices. 
Duke Math. J. 71 (1993), no. 1, 1--21. 

\bibitem{16} C.E. Kenig, G. Ponce and L. Vega,
A bilinear estimate with applications to the KdV equation. 
J. Amer. Math. Soc. 9 (1996), no. 2, 573--603.

\bibitem{9}
K. Nishikawa, H. Hojo, K. Mima and H. Ikezi,
Coupled nonlinear electron-plasma and ion acoustic waves,
Phys. Rev. Lett. 33 (1974), 148--151.


\bibitem{11} M. Tsutsumi,
Well-posedness of the Cauchy problem for a coupled Schršdinger-KdV equation. Nonlinear mathematical problems in industry, II (Iwaki, 1992), 513--528, 
GAKUTO Internat. Ser. Math. Sci. Appl., 2 (1993).

\bibitem{10}
N. Yajima and J. Satsuma, Soliton solutions in a diatomic lattice system,
Prog. Theor. Phys. 62 (1979), 370--378.



\end{thebibliography}


\end{document}